\documentclass{amsart}
\usepackage{amssymb,amsthm,amsmath,epsfig,latexsym,xypic}
\usepackage{calc}
\usepackage{latexsym}

\begin{document}

\newcommand{\mmbox}[1]{\mbox{${#1}$}}
\newcommand{\proj}[1]{\mmbox{{\mathbb P}^{#1}}}
\newcommand{\affine}[1]{\mmbox{{\mathbb A}^{#1}}}
\newcommand{\Ann}[1]{\mmbox{{\rm Ann}({#1})}}
\newcommand{\caps}[3]{\mmbox{{#1}_{#2} \cap \ldots \cap {#1}_{#3}}}
\newcommand{\N}{{\mathbb N}}
\newcommand{\Z}{{\mathbb Z}}
\newcommand{\R}{{\mathbb R}}
\newcommand{\Tor}{\mathop{\rm Tor}\nolimits}
\newcommand{\Ext}{\mathop{\rm Ext}\nolimits}
\newcommand{\Hom}{\mathop{\rm Hom}\nolimits}
\newcommand{\im}{\mathop{\rm Im}\nolimits}
\newcommand{\rank}{\mathop{\rm rank}\nolimits}
\newcommand{\codim}{\mathop{\rm codim}\nolimits}
\newcommand{\supp}{\mathop{\rm supp}\nolimits}
\newcommand{\CB}{Cayley-Bacharach}
\newcommand{\HF}{\mathrm{HF}}
\newcommand{\HP}{\mathrm{HP}}
\newcommand{\coker}{\mathop{\rm coker}\nolimits}
\sloppy
\newtheorem{defn0}{Definition}[section]
\newtheorem{prop0}[defn0]{Proposition}
\newtheorem{conj0}[defn0]{Conjecture}
\newtheorem{thm0}[defn0]{Theorem}
\newtheorem{lem0}[defn0]{Lemma}
\newtheorem{corollary0}[defn0]{Corollary}
\newtheorem{example0}[defn0]{Example}

\newenvironment{defn}{\begin{defn0}}{\end{defn0}}
\newenvironment{prop}{\begin{prop0}}{\end{prop0}}
\newenvironment{conj}{\begin{conj0}}{\end{conj0}}
\newenvironment{thm}{\begin{thm0}}{\end{thm0}}
\newenvironment{lem}{\begin{lem0}}{\end{lem0}}
\newenvironment{cor}{\begin{corollary0}}{\end{corollary0}}
\newenvironment{exm}{\begin{example0}\rm}{\end{example0}}

\newcommand{\aCM}{arithmetically Cohen-Macaulay }
\newcommand{\SES}{short exact sequence }
\newcommand{\LES}{long exact sequence }
\newcommand{\SRring}{Stanley-Reisner ring }
\newcommand{\SRrings}{Stanley-Reisner rings }

\newcommand{\msp}{\renewcommand{\arraystretch}{.5}}
\newcommand{\rsp}{\renewcommand{\arraystretch}{1}}
 
\newenvironment{lmatrix}{\renewcommand{\arraystretch}{.5}\small
 \begin{pmatrix}} {\end{pmatrix}\renewcommand{\arraystretch}{1}}
\newenvironment{llmatrix}{\renewcommand{\arraystretch}{.5}\scriptsize
 \begin{pmatrix}} {\end{pmatrix}\renewcommand{\arraystretch}{1}}
\newenvironment{larray}{\renewcommand{\arraystretch}{.5}\begin{array}}
 {\end{array}\renewcommand{\arraystretch}{1}}
 
\def \a{{\mathrel{\smash-}}{\mathrel{\mkern-8mu}}
{\mathrel{\smash-}}{\mathrel{\mkern-8mu}}
{\mathrel{\smash-}}{\mathrel{\mkern-8mu}}}

\newcommand{\std}{Gr\"{o}bner}
\newcommand{\jq}{J_{Q}}

\title {Betti numbers and degree bounds for some linked zero-schemes}

\author{Leah Gold}
\address{Gold: Mathematics Department \\ Texas A\&M University \\
  College Station \\ TX 77843-3368 \\ USA}
\thanks{Gold is supported by an NSF-VIGRE postdoctoral fellowship}
\email{lgold@math.tamu.edu}

\author{Hal Schenck}
\thanks{Schenck is supported by NSF Grant DMS 03--11142 and 
NSA Grant MDA 904-03-1-0006.}
\address{Schenck: Mathematics Department \\ Texas A\&M University \\
  College Station \\ TX 77843-3368\\ USA}
\email{schenck@math.tamu.edu}

\author{Hema Srinivasan}
\address{Srinivasan: Mathematics  Department \\ University of Missouri  \\ 
Columbia \\ MO 65211 \\ USA}
\email{hema@math.missouri.edu}
 

\begin{abstract}
\noindent
In \cite{HS}, Herzog and Srinivasan study the relationship 
between the graded Betti numbers of a homogeneous 
ideal $I$ in a polynomial ring $R$ and the degree of $I$. For 
certain classes of ideals, they prove a bound on the degree in terms 
of the largest and smallest Betti numbers, generalizing results
of Huneke and Miller in \cite{HM}. The bound is conjectured to 
hold in general; we study this using linkage. 
If $R/I$ is Cohen-Macaulay, we may reduce to the case where 
$I$ defines a zero-dimensional subscheme $Y.$ If $Y$ is residual
to a zero-scheme $Z$ of a certain type (low degree
or points in special position), then we show that the conjecture is
true for $I_Y$. 
\end{abstract}
\maketitle


\section{Introduction}\label{sec:one}
Let $R$ be a polynomial ring over a field $\mathbb{K}$, and let $I$ be a
homogeneous ideal. Then the module $R/I$ admits a finite minimal graded
free resolution over $R$:
\begin{center}
  $ \mathbb{F}: \mbox{ }\cdots \rightarrow \bigoplus\limits_{j \in
    J_2} R(-d_{2,j}) \rightarrow \bigoplus\limits_{j \in J_1}
  R(-d_{1,j}) \rightarrow R \rightarrow R/I\rightarrow 0.$
\end{center}
Many important numerical invariants of $I$ and the associated scheme can be
read off from the free resolution. For example, the {\em Hilbert
polynomial} is the polynomial $f(t) \in \mathbb{Q}[t]$ such that for all $m
\gg 0$, $\dim_\mathbb{K} (R/I)_m = f(m)$; if $f(t)$ has degree $n$ and lead
coefficient $d$, then the {\it degree} of $I$ is $n!d$.  When one has an
explicit free resolution in hand, then it is possible to write down the
Hilbert polynomial, and hence the degree, in terms of the shifts $d_{i,j}$
which appear in the free resolution.

If $R/I$ is Cohen-Macaulay and has a {\em pure resolution} 
\[ 0 \rightarrow R^{e_p}(-d_p)\cdots
  \rightarrow R^{e_2}(-d_2) \rightarrow R^{e_1}(-d_1)\rightarrow R
  \rightarrow R/I\rightarrow  0,
\]
then Huneke and Miller show in \cite{HM} that $deg(I) = (\prod_{i=1}^pd_i)/p!$.
Their result points to a more general possibility:

\begin{conj}[Huneke \& Srinivasan] \label{conj1}
  Let $R/I$ be a Cohen-Macaulay algebra with minimal free resolution of the
  form
  \begin{center}
  $ 0 \rightarrow \bigoplus\limits_{j \in J_p} R(-d_{p,j}) \rightarrow
  \dots \rightarrow \bigoplus\limits_{j \in J_2} R(-d_{2,j})
  \rightarrow \bigoplus\limits_{j \in J_1} R(-d_{1,j}) \rightarrow R 
  \rightarrow R/I \rightarrow 0$.
  \end{center}
  Let $m_i = \min \; \{d_{i,j} \;|\; j \in J_i\}$ be the minimum
  degree shift at the $i$th step and let $M_i = \max \;\{d_{i,j} \;|\:
  j \in J_i\}$ be the maximum degree shift at the $i$th step. Then
  \[
  \frac{\prod_{i=1}^p m_i}{p!}  \;\leq \; \deg(I) \;\leq\;
  \frac{\prod_{i=1}^p M_i}{p!}.
  \]
\end{conj}

When $R/I$ is not Cohen-Macaulay, it is easy to see that the lower
bound fails; for example if $I = (x^2, x y) \subset k[x,y]$, then
$\deg(I)= 1$, $m_1 = 2$ and $m_2 = 3$, but $\frac{(2)(3)}{2!} \geq 1$.
However, in \cite{HS}, Herzog and Srinivasan conjecture that even if
$R/I$ is not Cohen-Macaulay, the upper bound is still valid if one
takes $p = \codim(I)$. Conjecture 1.1 is verified in \cite{HS} in a
number of situations: when $I$ is codimension two; for codimension
three Gorenstein ideals with five generators (in fact, the upper bound
holds for codimension three Gorenstein with no restriction on the
number of generators); when $I$ is a complete intersection, and also
for certain classes of monomial ideals. Additional cases where
Conjecture 1.1 has been verified appear in \cite{GHP}, \cite{G},
\cite{GV}.  In the non-Cohen-Macaulay case, \cite{HS}
proves the bound for stable monomial ideals \cite{EK}, squarefree
strongly stable monomial ideals \cite{AHH}, and ideals with a pure
resolution; \cite{R} proves it for codimension two. In fact, in
the codimension two Cohen-Macaulay and codimension three Gorenstein
cases, a stronger version of the conjecture holds, see \cite{MNR}.

Most of the situations where the conjecture is known to be
true are when the entire minimal free resolution is known; the work in
proving the conjecture generally involves a complicated analysis
translating the numbers $d_{i,j}$ to the actual degree.  In this paper
we take a different approach. Our goal is to obtain {\it only}
the information germane to the conjecture; in particular we need the
smallest and biggest shift at each step. When $I$ is Cohen-Macaulay we
can always slice with hyperplanes without changing the degree or free
resolution, hence the study of the conjecture, in the Cohen-Macaulay
case, always reduces to the study of zero-schemes.

Suppose $Y$ is a zero-scheme, and $Z$ is a zero-scheme residual to $Y$
inside a complete intersection $X$. The resolution for $I_X$ is known, so
if one has some control over $Z$, (for example, when $Z$ consists of a
small number of points, or points in special position), then linkage allows
us to say something about the resolution for $I_Y$. 
Central to this are the results of Peskine-Szpiro \cite{PS} connecting 
resolutions and linkage.

\subsection{Resolutions and linkage}
Two codimension $r$ subschemes $Y$ and $Z$ of $\mathbb{P}^n$ are 
{\em linked} in a complete intersection $X$ if 
$I_Y = I_X:I_Z$ and $I_Z = I_X:I_Y$.
The most familiar form of linkage is the Cayley-Bacharach 
theorem \cite{EGH1}, which was our original motivation.

\begin{thm}[see \cite{PS} or \cite{N}]\label{linkres}
  Let $X \subset \proj{n}$ be an arithmetically
  Gorenstein scheme of codimension $n$, with minimal free resolution
\[
0 \rightarrow R(-\alpha) \rightarrow F_{n-1} \rightarrow F_{n-2}
\rightarrow \cdots \rightarrow F_1 \rightarrow R \rightarrow R/I_X
\rightarrow 0.
\]
Suppose that $Z$ and $Y$ are linked in $X$, and that the minimal 
free resolution of $R/I_Z$ is given by:
\[
0  \rightarrow  G_n \rightarrow G_{n-1} \rightarrow  \cdots \rightarrow
G_1 \rightarrow  R \rightarrow R/I_Z \rightarrow 0.
\]
 Then there is a free resolution for $R/I_Y$ given by 
\[
0  \rightarrow  G_1^\vee(-\alpha)  \rightarrow 
\begin{array}{c}
G_2^\vee(-\alpha) \\
\oplus\\ 
F_1^\vee(-\alpha)
\end{array}
 \rightarrow  
\begin{array}{c}
G_3^\vee(-\alpha) \\
\oplus\\ 
F_2^\vee(-\alpha)
\end{array}
 \rightarrow  
\cdots \rightarrow
\begin{array}{c}
G_n^\vee(-\alpha) \\
\oplus\\ 
F_{n-1}^\vee(-\alpha)
\end{array}
 \rightarrow  
R \rightarrow  
R/I_Y \rightarrow  0.
\]
\end{thm}

It turns out that in certain situations 
the shifts in the mapping cone resolution for $Y$ given by the
theorem above are such that no cancellation
of the relevant shifts can occur.

\section{Ideals linked to a collinear subscheme}
We assume for the remainder of the paper that $n \geq 3$ and that
$X$ is a non-degenerate (all the $d_i > 1$) complete intersection 
zero-scheme of type $(d_1,d_2,\ldots, d_n)$; let $d_X$ denote the
degree of $X$, and $\alpha_X = \sum_{i=1}^n d_i$. 
Suppose $Z$ is a complete intersection subscheme of $X$, of type
$(e_1,\ldots, e_n)$; with $d_Z$ and $\alpha_Z$ as above.
A minimal free resolution for $R/I_X$ is given by 
$F_i = \wedge^i(\oplus_{j=1}^n R(-d_j))$, and a minimal free
resolution for $R/I_Z$ is given by 
$G_i = \wedge^i(\oplus_{j=1}^n R(-e_j))$. 
In this case it is easy to see that Theorem 1.2 implies that 
there exists $f$ of degree $a = \alpha_X-\alpha_Z$ such that 
$I_Y = I_X:I_Z = (I_X+f)$ and $I_Z = I_X:f$; in particular, $I_Y$ is
an almost complete intersection. Since $I_X \subseteq I_Z$, $R/I_X
\rightarrow R/I_Z$; the mapping cone of Theorem 1.2 comes from
a map of complexes which begins:
\begin{small}
\[
\xymatrix{
 \ar[r] &\wedge^2(\oplus_{i=1}^n R(-d_i)) \ar[r]& \oplus_{i=1}^n
 R(-d_i) \ar[r]\ar[d]^{\phi} & R \ar[r] \ar[d] &
 R/I_X \ar[d] \ar[r]& 0\\
\ar[r]&  \wedge^2(\oplus_{i=1}^n R(-e_i)) \ar[r] &\oplus_{i=1}^n R(-e_i)\ar[r] &
  R \ar[r] & R/I_Z \ar[r] &0\\
}
\]
\end{small}
The comparison map $\phi$ 
which makes the diagram commute is simply an expression of the
generators of $I_X$ in terms of the generators of $I_Z$ 
(e.g.\cite{E}, Exercise 21.23). If 
$I_X \subseteq \mathfrak{m}I_Z$ then $\phi$ has entries in
$\mathfrak{m}$; in the construction of Theorem 1.2 the map
$G_{n-1}^\vee \rightarrow F_{n-1}^\vee$ is the transpose of $\phi$.
Since the comparison maps further back in the
resolution are simply exterior powers of $\phi$, we have:
\begin{lem}\label{nocancel}
 If $I_X \subseteq \mathfrak{m}I_Z$, then the mapping cone resolution
is in fact a minimal free resolution for $I_Y$.
\end{lem}
So if $I_X \subseteq \mathfrak{m}I_Z$, then the 
minimal free resolution $H_\bullet$ for $R/I_Y$ has
$H_n = \oplus_{i=1}^n R(e_i - \alpha_X))$, and for 
$i \in \{1,\ldots, n-1\}$, 
\[
H_i = \wedge^{n-i}(\oplus_{i=1}^n R(d_i)) \bigoplus
\wedge^{n-i+1}(\oplus_{i=1}^n R(e_i))(-\alpha_X).
\]
If $I_X \not\subseteq \mathfrak{m}I_Z$, then $I_X$ and $I_Z$
share some minimal generators; in this case, there can be 
cancellation in the mapping cone resolution:

\begin{exm} 
  Let $I_X = \langle x^2,y^2, z^6 \rangle \subseteq k[x,y,z,w]$,
  and let $I_Z = \langle x,y, z^6\rangle$. Then we find that $I_Y =
  I_X+\langle xy \rangle$. In betti diagram notation the mapping cone 
resolution of $R/I_Y$ is:
\begin{small}
$$
\vbox{\offinterlineskip 
\halign{\strut\hfil# \ \vrule\quad&# \ &# \ &# \ &# \ &# \ &# \ 
&# \ &# \ &# \ &# \ &# \ &# \ &# \
\cr
degree&1&4&6&3\cr
\noalign {\hrule}
0&1 &--&--&--&\cr
1&--&3 &2 &1 &\cr
2&--&--&1 &-- &\cr
3&--&--&--&--&\cr
4&--&--&--&--&\cr
5&--&1 &--&--&\cr
6&--&--&3 &2 &\cr
\noalign{\bigskip}
\noalign{\smallskip}
}}
$$
\end{small}
This is not a minimal resolution; the $R(-4)$ summand can be pruned off. 
The degree of $I_Y$ is $18$. Checking, we obtain 
$\prod_{i=1}^3 m_i = 54$, $\prod_{i=1}^3 M_i = 432$, and indeed
$9 \le 18 \le 72$. Notice that the upper bound was not affected
when we pruned the resolution, and the value of $\prod_{i=1}^3 m_i$
increased after pruning. 
\end{exm}
\begin{exm}\label{onept}
Let $Z$ be a single point. For $Y,$ Lemma~\ref{nocancel} implies that 
$M_n=m_n=\alpha_X-1$, and for $i<n$, 
$M_i = \alpha_X-n+i-1$ and $m_i = \sum_{j=1}^i d_j$ (where $d_i \le
d_j$ if $i \le j$). We want to show that
\[
( \prod_{j=1}^{n-1} \sum_{i=1}^j d_i) (\sum_{i=1}^n d_i -1) \;
\leq \; n!(d_X-1) \; \leq \; \prod_{i=1}^n (\alpha_X -i).
\]
For the upper bound there are two cases. If
$d_1 < d_n$, then we have the following inequalities:
\begin{eqnarray*}
n d_1 & \leq & d_1 + d_2 + \cdots +d_{n-1}+ d_n -1 = \alpha_X-1\\ (n-1) d_2 &
\leq & (d_2 + \cdots + d_n) + (d_1 -2) = \alpha_X -2\\ \vdots & & \vdots \\ 2
d_{n-1} & \leq & (d_{n-1} + d_n) + (d_1 + d_2 + \cdots + d_{n-2} - (n-1)) =
\alpha_X -(n-1)\\ d_{n} & \leq & (d_n) + (d_1 + d_2 + \cdots + d_{n-1} - n) =
\alpha_X -n
\end{eqnarray*}
So it follows that $n!(d_X-1) \; \leq \; n! d_1 d_2 \cdots d_n \; \leq
\; \prod_{i=1}^n (\alpha_X-i).$ If $d_1 = d_n = \delta$, then
\begin{eqnarray*}
n \delta & \leq & n \delta = \alpha_X \\
(n-1) \delta & \leq & (n-1) \delta + (\delta -2) = \alpha_X - 2(1)
  \; \leq \; \alpha_X -2 \\
(n-2) \delta & \leq & (n-2) \delta + (2)(\delta -2) = \alpha_X - 2(2)
 \; \leq \; \alpha_X - 3\\
\vdots & & \vdots \\
2 \delta & \leq & 2 \delta + (n-2)(\delta -2) = \alpha_X - 2(n-2)
  \; \leq \;  \alpha_X - (n-1)\\
\delta & \leq & \delta + (n-1)(\delta -2) = \alpha_X - 2(n-1)
\end{eqnarray*}
So $n! (\delta^n-1) \; \leq \; n! \delta^n \; \leq \;
(\alpha_X) \left( \prod_{i=2}^{n-1} (\alpha_X-i) \right) (\alpha_X -2n+2).$
To finish the upper bound, we must verify that 
$\alpha_X (\alpha_X -2n+2) \; \leq \; (\alpha_X -1)(\alpha_X -n)$;
this follows since $n \geq 3$.
 
The lower bound is easier: it holds for a complete intersection, and by assumption
$d_j \geq 2$ for all $j$, so we have
\[
\prod_{j=1}^{n} \sum_{i=1}^j d_i \le n!d_X \;\;\mbox{  and  }\;\;
j+1 \; \leq \; 2 j \; \leq \; \sum_{i=1}^j d_i.
\]
Thus
\[
n! = \prod_{j=1}^{n-1} (j+1) \; \leq \; \prod_{j=1}^{n-1} 2j
\; \leq \; \prod_{j=1}^{n-1} \sum_{i=1}^j d_i.
\]
Combining these two inequalities yields the lower bound. 
\end{exm}
\begin{lem}\label{ineq}
If $X$ is a non-degenerate zero-dimensional complete intersection in 
$\mathbb{P}^n$, with $n\ge 3$, then $d_X \leq { \alpha_X-1 \choose
  n}$, i.e. $d_X n! \le (\alpha_X-1)(\alpha_X-2) \cdots (\alpha_X-n)$.
\end{lem}
\begin{proof}
The bounds in  Conjecture~\ref{conj1} hold for a $(d_1, d_2,\cdots,
d_n)$ complete intersection, so 
$d_X \, n! \leq \alpha_X (\sum_{i=2}^n d_i) (\sum_{i=3}^n d_i) \cdots
  d_n.$ If $d_1 < d_n$, then as in the first case of Example~\ref{onept},
$d_X\, n! \leq (\alpha_X-1) (\sum_{i=2}^n
d_i) (\sum_{i=3}^n d_i) \cdots d_n$.
Hence it suffices to show
\[
 \alpha_X(\sum_{i=2}^n d_i) (\sum_{i=3}^n d_i) \cdots (\sum_{i=n}^n d_i)
 \leq \prod_{j=1}^{n}(\alpha_X-j)
\]

\noindent{\it Case 1:} $d_1 > 2$. Then $(\sum_{i=2}^n d_i) \leq
 (\alpha_X-3)$ and $(\sum_{i=j}^n d_i) \leq (\alpha_X-j)$ for all $j\geq 3$.  So
since $\alpha_X(\alpha_X-3) \leq (\alpha_X-1)(\alpha_X-2)$, we obtain:
\[
\begin{array}{rcl}
 \alpha_X (\sum_{i=2}^n d_i) (\sum_{i=3}^n d_i) \cdots (\sum_{i=n}^n d_i)
 &\leq &\alpha_X(\alpha_X-3) (\alpha_X-3) (\alpha_X -4) \cdots (\alpha_X-n)\\
 &\leq &\prod_{j=1}^{n}(\alpha_X-j)
\end{array}
\]
 
\noindent{\it Case 2:} $d_1 = 2$. Then $(\sum_{i=3}^n d_i) \leq (\alpha_X-4)$
 and $(\sum_{i=j}^n d_i) \leq (\alpha_X-j)$ for all $j\geq 2$, so
\[
 \alpha_X (\sum_{i=2}^n d_i) (\sum_{i=3}^n d_i) \cdots (\sum_{i=n}^n d_i)
 \leq
\alpha_X(\alpha_X-2) (\alpha_X-4) (\alpha_X -4) \cdots(\alpha_X-n).
\]
Since $\alpha_X (\alpha_X-4) \leq (\alpha_X-1)(\alpha_X-3)$, we obtain
$ \alpha_X(\alpha_X-2) (\alpha_X-4) (\alpha_X -4) \cdots(\alpha_X-n)
 \leq \prod_{j=1}^{n}(\alpha_X-j).$
\end{proof}
The proof of the next lemma is similar so we omit it.
\begin{lem}\label{betterineq}
With the same hypothesis as Lemma~\ref{ineq}, 
$d_X n! \leq \alpha_X (\alpha_X-2) (\alpha_X-4) (\alpha_X -6) \cdots
  (\alpha_X-2(n-1))$.
\end{lem}
\begin{defn}
  A subscheme $Z \subseteq \mathbb{P}^n$ is {\em collinear} if $I_Z =
  \langle l_1,\ldots, l_{n-1},f \rangle$, where the $l_i$ are linearly
  independent linear forms
  and $\deg f = t$.
\end{defn}
We now use linkage to study the case where $Y$ is linked in $X$ to a
collinear subscheme $Z$. While we expect our methods to work more
generally, this case is already complicated enough to be interesting.
Since the line $V(l_1,\ldots,l_{n-1})$ cannot be contained 
in each of the hypersurfaces defining $X$ (or $X$ would contain 
the whole line), the line on which $Z$
is supported must intersect one of the hypersurfaces defining
$X$ in a zero-scheme. Thus, $Z$ is of degree at most $d_n$.
Henceforth we write $\alpha$ for $\alpha_X$. 
\begin{thm}\label{ptsonlinethm}
  Let $X$ be a zero-dimensional complete intersection of type $d_1,
  d_2, \ldots, d_n$ in $\proj{n}$. Let $Z \subset X$ be a collinear
  subscheme of degree $t$, and let $Y$ be residual to $Z$. Then
  Conjecture~\ref{conj1} holds for $R/I_Y$.
\end{thm}
\begin{proof}
{\bf Upper bound.} 
  Because $d_j \geq 2$ for all $j$, even if cancellation occurs we
  have $M_i = \alpha -n+i-1$ for $i \in \{2, \ldots, n\}$, as in
  Example~\ref{onept}. For $i=1$, $M_1 \geq d_n$ or $M_1=d_n-1$, depending on
  the amount of cancellation. If $t \leq \sum_{i=1}^{n-1} (d_i-1)$,
  then $\alpha-n-t+1 \geq d_n$ and so $M_1 \geq d_n$. If
  $\sum_{i=1}^{n-1} (d_i-1) < t$, then cancellation can occur.

\smallskip
\noindent{\it Case 1: $M_1 \geq d_n$.} In this case, since
\[
n!(d-t) \leq n!d \leq \alpha (\alpha-d_1) (\alpha-d_1-d_2) \cdots (d_n),
\]
it suffices to show that 
\[
\alpha (\alpha-d_1) (\alpha-d_1-d_2) \cdots (d_{n-1}+d_n)(d_n)
\leq
(\alpha-1)(\alpha-2)(\alpha-3) \cdots (\alpha-(n-1))M_1
\]
Since $d_j \geq 2$ for all $j$, $\alpha(\alpha-d_1-d_2)
\leq (\alpha-1)(\alpha-3)$, and  
\[
\begin{array}{rcl}
(\alpha-d_1) &\leq &(\alpha-2)\\
(\alpha-d_1-d_2-d_3)  &\leq & (\alpha-4)\\
(\alpha-d_1-d_2-d_3-d_4) &\leq &(\alpha-5)\\
                         & \vdots, &
\end{array}
\] 
the result follows if $n\ge 5$. If $n = 4$, then we must replace
the $\alpha-4$ above with $M_1$. The result holds since $M_1 \ge d_4 =
\alpha-d_1-d_2-d_3$.

For $n=3$, there are four cases to analyze.
If $d_1 \geq 3$, then $\alpha(\alpha-d_1) \leq (\alpha-1)(\alpha-2)$.  
If $d_1 = 2$, then if $d_2 \geq 3$ we find that $6d \leq
(\alpha-1)(\alpha-2)d_3$ because $11 d_2 \leq d_2^2 + 2d_2d_3+d_3^2 +d_3$.
If $d_1=2$ and $d_2=2,$ but $d_3 \geq 3$, then we find that $24 d_3 \leq
d_3^3+5d_3^2 + 6d_3$. Since $d_3 \geq 3$, $18 \leq d_3^2+5d_3$ so the
inequality is true.

Finally, if $d_1=d_2=d_3=2$, then as long as $t>1$ we have $6(8-t) \leq
(5)(4)(2)$, so the bound holds when $t>1$. The case $t=1$ is covered
by Example~\ref{onept}, which concludes Case 1.

\smallskip
\noindent{\it Case 2: $d_n > M_1$.} Then $\alpha-t-n+1 = d_n-1$.  If $d_1 =
d_n$, then since at most $n-1$ of the $d_i$'s can cancel, this forces $M_1
= d_1 =d_n$ and the inequalities from the previous case apply. So
henceforth we assume $d_1 < d_n$, which as noted in Lemma~\ref{ineq}
implies $d \, n!  \leq (\alpha-1) (\sum_{i=2}^n d_i) (\sum_{i=3}^n d_i)
\cdots d_n$. We wish to show
\[
n!(d-t) \leq (\alpha-t-n+1) \prod_{i=2}^n (\alpha-n+i-1)
= (\alpha-t-n+1) \prod_{i=1}^{n-1} (\alpha-i)
\]
Suppose $n \geq 5$. We claim that $d_n(d_n+d_{n-1}) \leq (d_n-1)(\alpha - n+2)
= (d_n-1)(d_n+t)$.  This follows from the inequalities
\[
\begin{array}{rcl}
(d_n-1)(d_n+t) - d_n(d_n+d_{n-1}) &=& -d_n + t(d_n-1) - d_{n-1}d_n\\ &\geq
&-d_n + (d_n-1)(d_{n-1} + n-2) -d_{n-1}d_n
\end{array}
\]
because $t = \alpha-d_n -n +2 = d_{n-1} + \sum_{i=1}^{n-2} (d_i-1) \geq
d_{n-1}+n-2$. Then
\[
\begin{array}{rcl}
-d_n + (d_n-1)(d_{n-1} + n-2) -d_{n-1}d_n & = &-d_n + (n-2) d_n -
 d_{n-1}-(n-2)\\
 &=& (n-4)d_n + (d_n-d_{n-1}) - (n-2)\\
 &\geq& (n-4) d_n - (n-2)\\
 &=& (n-4) (d_n-1) -2.
\end{array}
\]
Finally $(n-4)(d_n-1) \geq 2$ because $n \geq 5$ and $d_n > d_1 \geq 2$, so
we obtain
\begin{small}
\[
\begin{array}{rcl}
n!d &\leq &d_n (d_n+d_{n-1}) (d_n+d_{n-1}+d_{n-2})  \cdots
(\alpha-d_1)(\alpha-1)\\
& \leq &(d_n-1) (d_n+t) (d_n+d_{n-1}+d_{n-2})  \cdots
(\alpha-d_1)(\alpha-1)\\
&= &(d_n-1) (\alpha-n+2) (d_n+d_{n-1}+d_{n-2}) \cdots
(\alpha-d_1)(\alpha-1)\\ 
&\leq& (d_n-1) (\alpha-n+2) (\alpha-n+1)(d_n+d_{n-1}+d_{n-2}+d_{n-3}) \cdots
(\alpha-d_1)(\alpha-1)\\ 
&\leq &(d_n-1) (\alpha-n+2) (\alpha - n+1) 
(\alpha-(n-3))(\alpha-(n-4))  \cdots (\alpha-2)(\alpha-1).
\end{array}
\]
\end{small}
Hence, the upper bound holds if $n\ge 5$. 

\smallskip
If $n= 4$ and $d_2 < d_4$, then $3d_2 \leq d_2+d_3+d_4-1+d_1-2 = \alpha
-3$. If $d_4=d_3$, then since $d_1<d_4$, we also have $4d_1 \leq \alpha-2$.
So, $12d_1d_2 \leq (\alpha -2)(\alpha -3)$. On the other hand, if
$d_2=d_4$, then $3d_2 \leq \alpha-2$ and $4d_1 \leq \alpha-3$ so we also
find that $12d_1d_2 \leq (\alpha -2)(\alpha -3)$. It just remains to show
that $2d_3d_4 \leq (\alpha -1)(d_4-1)$.  But $(d_4-1)(\alpha -1) -2d_3d_4
\ge (d_4-1)(2d_4+3)-2d_4^2=d_4-3\ge 0$.  Thus the upper bound holds when
$d_4=d_3$.  If $d_3<d_4$, we may only have $4d_1\leq (\alpha -1)$.
Nevertheless,
\[
\begin{array}{rcl}
(\alpha -2)(d_4-1)-2d_3d_4 &= &(d_1+d_2+d_4-d_3-2)(d_4-1)-2d_3 \\
& \ge & (d_1+d_2+d_4-d_3- 2)(d_4-1)-2(d_4-1)\\
& = &(d_1+d_2+d_4-d_3-4)(d_4-1)\\
& = & (d_1+d_2-4 +d_4- d_3)(d_4-1) \ge 0.
\end{array}
\]
Thus, the upper bound holds when $n=4$.

\smallskip
If $n=3$, then since $M_1 =d_3-1$, $d_2 \neq d_3$.  If $3d_1\leq (\alpha
-2)$ then as before, $(\alpha -1)(d_3-1)-2d_2d_3 \ge (d_1-d_2+d_3-3)(d_3-1)
\ge 0$.  If $3d_1=\alpha -1$, we must have $d_1=d_2=d_3-1$.  In this case,
using the fact that $t=2d_1-1$, we calculate the inequality directly:
$6(d_1^2(d_1+1)- (2d_1-1) ) \leq (d_1)(3d_1+1-2) (3d_1+1-1)$ simplifies to
the true statement $0 \leq 3(d_1-1)(d_1-2d_1+2)$.
\vskip .15in

\noindent {\bf Lower bound.}
If there is no cancellation, then $m_n = \alpha - t$ and for $i<n$ 
we have $m_i = \min\{\alpha-n-t+i, \sum_{j=1}^i d_j \}$.
In particular, $m_i \leq \sum_{j=1}^i d_j$, for $i \in \{1,\ldots,
n-1 \}$, and so
\[
\prod_{i=1}^n m_i = (\prod_{i=1}^{n-1} m_i) m_n \leq (\prod_{i=1}^{n-1}
\sum_{j=1}^i d_j) m_n.
\]
Hence it is sufficient to prove that
\[
(\prod_{i=1}^{n-1} \sum_{j=1}^i d_j) (\alpha-t) \leq n!(d-t) .
\]
Exactly as in Example~\ref{onept}, we have 
\[
\prod_{i=1}^{n-1} \sum_{j=1}^i d_j \alpha \le n!d \;\;\mbox{  and  }\;\;
i+1 \; \leq \; 2 i \; \leq \; \sum_{j=1}^i d_j.
\]
So $n! t  = t \prod_{i=1}^{n-1} (i+1) \; 
\leq \; t \prod_{i=1}^{n-1} \sum_{j=1}^i d_j.$
Subtracting this inequality from the left hand inequality above yields
the desired inequality, so the lower bound holds for $R/I_Y$ if
there is no cancellation.

\smallskip
Now let us look at where cancellation can occur. We only care about
cancellation when a term of some degree that shows up in the set of
minimums disappears. We can break it up into two cases:

\smallskip 
\noindent{\it Case 1: $t < d_n$.} Then $\alpha-t > \alpha-d_n$, and so
$\alpha-t -1 \ge \alpha-d_n$, hence $m_{n-1} \le \alpha-d_n$. Also
$\alpha-t -1 \ge \alpha-d_n$ implies $\alpha-t -1 >
\alpha-d_n-d_{n-1}$, so that $m_{n-2} \le \alpha-d_n-d_{n-1}$, and in
general $m_{n-i} \le \alpha-d_n-\cdots d_{n-i+1}$. So if $m_n =
\alpha-t$, then the argument from the previous case holds.  

However, if $t=d_l$ for some $l<n$, then it is possible that $m_n =
\alpha-1$.  So in this case, we need to show that
\[
(\prod_{i=1}^{n-1} \sum_{j=1}^i d_j) (\alpha-1)
\leq n! (d-d_l). 
\]
We have the inequalities
\[
\begin{array}{ccc}
d_1 & \le & d_1 \\
d_1+d_2 & \le & 2d_2 \\ 
\vdots\\ 
d_1+d_2+\cdots +d_{n-2}+d_{n-1} & \le & (n-1)d_{n-1} \\ 
\alpha+1 & \le & nd_n,
\end{array}
\]
where the last row follows since $d_l < d_n$. 
Subtracting $2\prod_{i=1}^{n-1} \sum_{j=1}^i d_j$ from the
product of the left hand column and $n!d_l$ from the product of
the right hand column would yield the desired inequality, so
it suffices to show that $n!d_l \le 2\prod_{i=1}^{n-1} \sum_{j=1}^i
d_j$. Let $\beta = \prod_{i=1}^{n-2} \sum_{j=1}^i d_j$, so 
\[
2\prod_{i=1}^{n-1} \sum_{j=1}^i d_j  = 2(d_{n-1}+\sum_{j=1}^{n-2} d_j)\beta. 
\]
Since $d_l \le d_{n-1}$, it is enough to show that $n! \le
2\beta$. Since the $d_i$ are at least two, 
\[
2^{n-1}(n-2)! \le 2\beta, 
\]
and the inequality holds if $n\ge 6$. For $n \in \{3,4,5\}$, a 
case analysis shows we have to verify the bound directly for
\[
\begin{array}{ll}
n=3 & d_1=2 \\
n=4 & (d_1,d_2) = (2,2) \mbox{ or }(2,3) \\ 
n=5 & (d_1,d_2,d_3) = (2,2,2) \mbox{ or }(2,2,3).
\end{array}
\]
For example, if $n=3$ and $d_1=2$, we must verify that 
\[
2(2+d_2)(2+d_2+d_3-1) \le 6(2d_2d_3-d_2). 
\]
This follows by summing the inequalities:
\[
\begin{array}{ccc}
(2+d_2)d_3 & \le & (2d_2)d_3 \\
(2+d_2)(d_2+1) & \le & (2d_2)d_3,
\end{array}
\]
and observing that $2d_2d_3-3d_2=d_2(2d_3-3)\ge 0$. The other
cases are similar so we omit them.

\smallskip 
\noindent{\it Case 2: $t = d_n$.} The $\alpha-d_n$ term cancels with $\alpha-t$, and
so $m_n = \alpha-1$. Also $m_{n-1} = \min\{\alpha-d_{n-1}, \alpha-t-1\}
\leq \alpha-t -1 = \alpha - d_n-1$. Since all the $d_i \geq 2$, we cannot
have $\alpha-d_n- \cdots -d_{k+1} = \alpha-n-t+k+1$ for any $k \leq n-2$,
and hence we always have $m_i \leq \sum_{j=1}^i d_j$ for $i \leq n-2$. In
order to prove the lower bound, we need to show
\[
(\alpha-1)(\alpha-d_n-1) \prod_{i=1}^{n-2} \sum_{j=1}^i d_j
\leq
n! (d-d_n) 
\]
We can write 
\[
n! (d-d_n) = d_n n (n-1)! (d'-1)
\]
where $d' = \prod_{i=1}^{n-1} d_i$. By the bound on the complete
intersection of type $d_1, d_2, \ldots, d_{n-1}$, we know that
\[
(n-1)! d' \geq \prod_{i=1}^{n-1} \sum_{j=1}^i d_j = (\alpha - d_n)
\prod_{i=1}^{n-2} \sum_{j=1}^i d_j.
\]
It is also true that $n-1 \leq 2^{n-2}$ for all $n \geq 2$, so 
\[
(n-1)! \leq 2^{n-2} (n-2)! \leq \prod_{i=1}^{n-2} \sum_{j=1}^i d_j,
\mbox{ since }d_i \ge 2.
\]
Therefore
\[
(n-1)! (d'-1) = (n-1)! d' - (n-1)! \geq (\alpha-d_n) \prod_{i=1}^{n-2}
\sum_{j=1}^i d_j - \prod_{i=1}^{n-2} \sum_{j=1}^i d_j = (\alpha-d_n-1)
\prod_{i=1}^{n-2} \sum_{j=1}^i d_j.
\]
But since $nd_n \geq \alpha \ge \alpha-1$, this gives
\[
n! (d - d_n) = d_n n (n-1)! (d'-1) \geq (\alpha-1) (\alpha-d_n-1)
\prod_{i=1}^{n-2} \sum_{j=1}^i d_j.
\]
\end{proof}

\section{$Y$ is linked to 3 general points}
In this section, we study the simplest $Z$ which is not a collinear scheme:
three general points. While we are able to carry out the degree analysis in
this case, it also serves to illustrate that this type of argument
will become increasingly complex. 

\begin{thm}\label{3noncollinearptsthm}
  Let $X$ be a zero-dimensional complete intersection of type $d_1, d_2,
  \ldots, d_n$ in $\proj{n}$, $n>2$. Let $Z \subset X$ be a set of 3
  non-collinear points, and suppose $Y$ is linked to $Z$ in $X$.  Then
  Conjecture~\ref{conj1} holds for $R/I_Y$.
\end{thm}
By Theorem~\ref{linkres}, the mapping cone resolution of $I_Y =
I_X:I_Z$ is
\begin{tiny}
\[
 0 \rightarrow 
 \begin{array}{c}
  R^{n-2} (-(\alpha-1))\\
  \oplus\\
  R^3 (-(\alpha-2))
 \end{array}
 \rightarrow 
 \cdots \rightarrow 
 \begin{array}{c}
  R^{{{n-2} \choose {n-i+1}}} (-(\alpha-n-1+i))\\
  \oplus\\
  R^{3 {{n-2} \choose {n-i}} + 2{{n-2} \choose {n-i-1}}}  (-(\alpha-n-2+i))\\
  \oplus\\
  \oplus R (-(\sum\limits_{j \in \Lambda \atop |\Lambda| = i} d_j))
 \end{array}
 \rightarrow 
 \cdots 
\]

\[
 \cdots \rightarrow 
 \begin{array}{c}
  R^{{n-2 \choose n-2}} (-(\alpha-n+2))\\
  \oplus\\
  R^{3 {n-2 \choose n-3} + 2{n-2 \choose n-4} } (-(\alpha-n+1))\\
  \oplus\\
  \oplus R (-(\sum\limits_{j \in \Lambda \atop |\Lambda| = 3} d_j))
 \end{array}
 \rightarrow 
  \begin{array}{c}
  R^{3 {n-2 \choose n-2} + 2{n-2 \choose n-3}}  (-(\alpha-n))\\
  \oplus\\
  \oplus R (-(d_j + d_k))
 \end{array}
 \rightarrow 
  \begin{array}{c}
  R^{ 2{n-2 \choose n-2}}  (-(\alpha-n-1))\\
  \oplus\\
  \oplus R (-d_j)
 \end{array}
 \rightarrow
 I_Y 
\]
\end{tiny}
\begin{proof}
{\bf Upper bound.} 
We begin with the upper bound. If $n\geq 4$, then there is no cancellation
of terms which affect the upper bound, and for $i \in \{3, \ldots, n-1 \}$, 
$M_i =  \max \{\sum_{j=n-i+1}^n d_j, \alpha-n+i-1\} = \alpha-n+i-1$, 
while 
\[
\begin{array}{l}
M_1 = \max \{ d_n, \alpha-n-1\} = \alpha-n-1\\
M_2 = \max \{d_{n-1}+d_n, \alpha-n\} = \alpha-n \\
M_n = \alpha - 1.
\end{array}
\]
So we want to show that
\[
n! (d-3) \leq (\alpha-n)(\alpha-(n+1)) \prod_{i=1}^{n-2} (\alpha-i).
\]
Since we know that 
\[
n!(d-3) \leq n! d \leq \alpha \prod_{i=2}^n \sum_{j=i}^n d_j,
\]
it is enough to show that 
\[
\alpha \prod_{i=2}^n \sum_{j=i}^n d_j \leq
(\alpha-n)(\alpha-(n+1)) \prod_{i=1}^{n-2} (\alpha-i).
\]

By Lemma~\ref{betterineq}, we know that
\[
\alpha \prod_{i=2}^n \sum_{j=i}^n d_j \leq
\alpha (\alpha-2) (\alpha-4) (\alpha -6) \cdots
(\alpha-2(n-1)),
\]
so it is enough to show that
\[
\alpha (\alpha-2) (\alpha-4) (\alpha -6) \cdots (\alpha-2(n-1))
\leq
(\alpha-n)(\alpha-(n+1)) \prod_{i=1}^{n-2} (\alpha-i).
\]

If $n>4$, then 
\[
\begin{array}{rcl}
\alpha-2 &\leq&\alpha-2\\
\alpha-2(3) &\leq& \alpha - 4\\
\alpha-2(4) &\leq& \alpha - 5\\
&\vdots &\\
\alpha-2(n-3) &\leq& \alpha - (n-2)\\
\alpha-2(n-2) &\leq& \alpha-n\\
\alpha-2(n-1) &\leq& \alpha-(n+1)
\end{array}
\]
and
\[
\begin{array}{rcl}
\alpha(\alpha-4) &\leq&(\alpha-1)(\alpha-3).
\end{array}
\]
Taking the product, we see that the bound holds if $n>4$. If $n = 4$, 
then we must show that 
\[ 
\alpha(\alpha-2)(\alpha-4)(\alpha-6)
\leq (\alpha-1)(\alpha-2)(\alpha-4)(\alpha-5);
\] 
which is true since $\alpha(\alpha-6) \leq (\alpha-1)(\alpha-5)$ for all
$\alpha$.

Finally, if $n=3$, then we have to be a bit more careful. It is always true
that $M_1=\alpha-4$ and $M_3 = \alpha-1$. The value of $M_2$ is either
$\alpha-2$ or $\alpha-3$ depending on cancellation.

\noindent{\it Case 1: $d_1=d_2=d_3=2$.}
We check directly that
\[
30=3!(8-3) = (2)(6-3)(5) =(\alpha-4)(\alpha-3)(\alpha-1)\leq M_1M_2M_3.
\]
{\it Case 2: $d_1=d_2=2, d_3>2$.} 
In this case $\alpha=d_3+4$, and so $M_2 \geq d_3+1$. Again we plug in
values, and check to see that the resulting inequality is true. Is
$6(d-3) = 6(4d_3-3) \leq (d_3)(d_3+1)(d_3+3)$?  This is equivalent to
$0 \leq d_3^3+4d_3^2 -21 d_3 + 18 =(d_3-2)(d_3^2+6d_3-9)$, which is
true for $d_3 \geq 3$.

\noindent{\it Case 3: $d_1=2, d_2>2$.} Here $\alpha=d_2+d_3+2$ and $M_2 \geq
d_2+d_3-1$, so we need to check that $6(2d_2d_3-3) \leq
(d_2+d_3-2)(d_2+d_3-1)(d_2+d_3+1)$. This inequality reduces to checking
that $d_2^3+3d_2^2d_3+3d_2d_3^2+d_3^3-2d_2^2-16d_2d_3-2d_3^2-d_2-d_3+20
\geq 0$, which is true since for $3 \leq d_2 \leq d_3$,
\[
\begin{array}{rcl}
d_2^3+3d_2^2d_3+3d_2d_3^2+d_3^3 &\geq & 3d_2^2+9d_2d_3 +9d_3^2+3d_3^2\\
&=&2d_2^2+2d_3^2+d_2^2+9d_2d_3+8d_3^2\\
&\geq&2d_2^2+2d_3^2+d_2^2+9d_2d_3+7d_2d_3+d_3^2\\
&=&2d_2^2+2d_3^2+16d_2d_3+d_2^2+d_3^2\\
&\geq &2d_2^2+2d_3^2+16d_2d_3+d_2+d_3.
\end{array}
\]
{\it Case 4: $d_1>2$.}  In this case, we check directly that
\[
\alpha(\alpha-d_1)(\alpha-d_1-d_2) \leq
(\alpha-1)(\alpha-d_1)(\alpha-4) \leq M_1 M_2 M_3.
\]
The left expression is the familiar product from $I_X$, so it is bigger
than $3!d$, and hence also $3!(d-3)$. So the upper bound holds.
\vskip .1in
\noindent{\bf Lower bound}
Now we will prove the lower bound. Notice that the only cancellation
that is numerically feasible is at the last step because $d_j \geq 2$
for all $j$. So cancellation can only happen if $d_1$, $d_2$, and
possibly $d_3$ are all 2. Such a cancellation will affect $m_n$ only
if all three terms of degree $\alpha-2$ cancel, that is, if
$d_1=d_2=d_3=2$ and all possible cancellation occurs, and $d_4 \ge 3$
when $n\ge 4$. Therefore for $i< n$ we have 
$m_i=   \min \{\sum_{j=1}^i d_j,\alpha-n+i-2\}$, and $m_n$ is
either $\alpha - 1$ or $\alpha - 2$. 
If we assume $m_n=\alpha-2$, then there are four cases to consider.

\smallskip
\noindent{\it Case $n \geq 4$:}
We know that 
\[
(\alpha-2) \prod_{i=1}^{n-1} m_i 
\leq
(\alpha-2) \prod_{i=1}^{n-1} \sum_{j=1}^i d_j,
\]
so we need to show that the rightmost expression is less than or equal
to $n! (d-3)$. Since $d_j \geq 2$, $2i \leq \sum_{j=1}^i d_j$, so  
$n! 3 
\leq 
2^n (n-1)! \leq
2 \prod_{i=1}^{n-1} \sum_{j=1}^i d_j$. Thus
\[
(\alpha-2) \prod_{i=1}^{n-1} \sum_{j=1}^i d_j
\leq n!d - 2 \prod_{i=1}^{n-1} \sum_{j=1}^i d_j 
\leq n! d - n! 3 = n!(d-3) 
\]

\smallskip 
\noindent{\it Case $n=3$, $d_1=d_2=d_3=2$:} 
In this case $m_1 = 2$, $m_2 =3$, and $m_3=4$, so we check directly
that
\[
24 = (2)(3)(4) \leq 3! (2^3-3) = 30.
\]

\smallskip
\noindent{\it Case $n=3$, $d_1=d_2=2$, $d_3 > 2$:} 
In this case we check directly that $(2)(4)(\alpha-2) \leq 3!(d-3)$.
Since $\alpha= d_3+4$, this inequality holds as long as $d_3 \geq
\frac{17}{8}$, which it is.

\smallskip \noindent{\it Case $n=3$, $d_2 > 2$:} In this case $m_1 \leq d_1$,
$m_2 \leq d_1+d_2$, and $m_3=\alpha-2$.  Using the bound for the
complete intersection of type $d_1,d_2,d_3$, we have that
\[
d_1 (d_1+d_2)(\alpha-2) = 
d_1 (d_1+d_2)\alpha - 2 d_1 (d_1+d_2) \leq 3! d -18,
\]
which is true if $2 d_1 (d_1+d_2) \geq 18$. But $2 d_1 (d_1+d_2) \geq
2(2)(5) = 20$, so the bound holds.

If on the other hand $m_n = \alpha-1$, then it
must be true that $d_1=d_2=d_3=2$.  We know
\[
\prod_{i=1}^{n} m_i \leq (\alpha-1) \prod_{i=1}^{n-1} \sum_{j=1}^i d_j,
\]
and so it suffices to show
\[
\alpha \prod_{i=1}^{n-1} \sum_{j=1}^i d_j - \prod_{i=1}^{n-1} \sum_{j=1}^i
d_j \leq n!d - 3 \,n!,
\]
which would follow from
\[
3 \,n!
\leq \prod_{i=1}^{n-1} \sum_{j=1}^i d_j.
\]

Since $d_4 \geq 2$, we have that
\[
5!3 = 3 \cdot 2  \cdot 3 \cdot 4 \cdot 5
\leq
2 \cdot 4 \cdot 6 \cdot 8 
\leq
2 \cdot 4 \cdot 6 \cdot (6+d_4) = \prod_{i=1}^4 \sum_{j=1}^i d_j,
\]
and once $n$ is at least $6$, 
$\prod_{i=6}^n i \leq \prod_{i=5}^{n-1} \sum_{j=1}^i
d_j$; hence the desired inequality follows if $n \geq 5$.

If $n=4$, then we check directly. We have that $m_1=2$, $m_2 = 4$,
$m_3= 6$ and $m_4 = d_4+5$. A simple calculation shows
that in fact $4!(8d_4-3) \geq (2)(4)(6)(d_4+5)$ since $d_4 \geq 2$.

If $n=3$, then again we may check directly. We have that $m_1=2$,
$m_2=3$, and $m_3=5$. So we see that $30=3!(8-3) \geq (2)(3)(5)=30$.
\end{proof}
\vskip .15in

\noindent{\bf Acknowledgments} 
Macaulay~2 computations provided evidence for the results in this
paper. The first author thanks the University of Missouri
for supporting her visit during the fall of 2003, when portions of
this work were performed. 

\renewcommand{\baselinestretch}{1.0}
\small\normalsize 
 
\bibliographystyle{amsalpha}

\end{document}